\documentclass[11pt]{article}
\usepackage[latin1]{inputenc}
\usepackage{amsmath,amsfonts,amscd,amssymb}
\usepackage[all]{xy}
\usepackage{graphicx,epsfig}
\usepackage{color}
\usepackage[]{hyperref}
\usepackage{makeidx}
\newtheorem{thm}{Theorem}[section]
\newtheorem{cor}{Corollary}[section]

\newtheorem{prop}{Proposition}[section]
\newtheorem{deft}{Definition}[section]
\newtheorem{rek}{Remark}[section]

\let\noi=\noindent

\newcommand{\dsp}{\displaystyle}

\def\N{\mathbb{N}} 
\def\Z{\mathbb{Z}}  
\def\R{\mathbb{R}}  
\def\Q{\mathbb{Q}}
\def\F{\mathbb{F}}
\def\sper{\mbox{Sper}}

\def\notin{\mbox{$\in$ \hspace{-.8em}/}} 
\def\sgn{\mbox{sgn}} 

\title{On points at infinity of real spectra of polynomial rings}
\author{F. Lucas\\
D\'{e}partement de Math\'{e}matiques\\
Universit\'{e} d'Angers\\
2, bd Lavoisier\\
49045 Angers Cedex, France \and
D. Schaub\\
D\'{e}partement de Math\'{e}matiques\\
Universit\'{e} d'Angers\\
2, bd Lavoisier\\
49045 Angers Cedex, France \and
M. Spivakovsky\\
CNRS-Institut de Math\'{e}matiques de Toulouse\\
Universit\'e Paul Sabatier\\
118, route de Narbonne\\
31062 Toulouse Cedex 9, France.}

\date{}

\begin{document}

\maketitle

\begin{abstract} Let $R$ be a real closed field and $A=R[x_1,\dots,x_n]$. Let $\sper\ A$ denote the real spectrum of $A$.
There are two kinds of points in $\sper\ A$: finite points (those for which all of $|x_1|$,\dots,$|x_n|$ are bounded above
by some constant in $R$) and points at infinity. In this paper we study the structure of the set of points at infinity of
$\sper\ A$ and their associated valuations. Let $T$ be a subset of $\{1,\dots,n\}$. For $j\in\{1,\dots,n\}$, let $y_j=x_j$
if $j\notin T$ and $y_j=\frac1{x_j}$ if $j\in T$. Let $B_T=R[y_1,\dots,y_n]$. We construct a finite partition $\sper\
A=\coprod\limits_TU_T$ and a homeomorphism of each of the sets $U_T$ with a subspace of the space of finite
points of $\sper\ B_T$. For each point $\delta$ at infinity in $U_T$, we describe the associated valuation $\nu_{\delta^*}$
of its image $\delta^*\in \sper\ B_T$ in terms of the valuation $\nu_\delta$ associated to $\delta$. Among other things we
show that the valuation $\nu_{\delta^*}$ is composed with $\nu_\delta$ (in other words, the valuation ring $R_\delta$ is a
localization of $R_{\delta^*}$ at a suitable prime ideal).
\end{abstract}

\section{Introduction}

Let $R$ be a real closed field and $z_0,\dots,z_n$ independent variables. A basic fact of life in mathematics is the way
the $n$-dimensional projective space $\mbox{Proj}\ R[z_0,\dots,z_n]$ and other rational projective schemes such as
$\left(\mathbb P_R^1\right)^n$ are glued together from affine charts of the form $\mbox{Spec}\ R[x_1,\dots,x_n]$.
Given two such coordinate charts $\mbox{Spec}\ R[x_1,\dots,x_n]$ and $\mbox{Spec}\ R[y_1,\dots,y_n]$, it is often easy to
write down formulae describing the coordinate transformation from the $x$ to the $y$ coordinates. The subject of
this paper is a part of the analogous story for real spectra (see Definition \ref{sper} below), which is more interesting,
because the real spectrum $\sper\ R[x_1,\dots,x_n]$ already contains much information ``at infinity''.

To explain this in more detail, we first recall the definition of real spectrum and other related objects, studied in this
paper.\medskip

\noi\textbf{Notation and conventions.} All the rings in this paper will be commutative with 1. For a prime ideal $\mathfrak
p$ in a ring $B$, $\kappa(\mathfrak p)$ will denote the residue field of the local ring $B_{\mathfrak p}$:
$\kappa(\mathfrak p)=\frac{B_{\mathfrak p}}{\mathfrak pB_{\mathfrak p}}$.\medskip

Let $B$ be a ring. A point $\alpha$ in the real spectrum of $B$ is, by
definition, the data of a prime ideal $\mathfrak{p}$ of $B$, and a
total ordering $\leq$ of the quotient ring $B/\mathfrak{p}$, or,
equivalently, of the field of fractions of $B/\mathfrak{p}$.  Another
way of defining the point $\alpha$ is as a homomorphism from $B$ to a
real closed field, where two homomorphisms are identified if they have
the same kernel $\mathfrak{p}$ and induce the same total ordering on
$B/\mathfrak{p}$.

The ideal $\mathfrak{p}$ is called the support of $\alpha$ and denoted
by $\mathfrak{p}_\alpha$, the quotient ring $B/\mathfrak{p}_{\alpha}$
by $B[\alpha]$, its field of fractions by $B(\alpha)$ and the real
closure of $B(\alpha)$ by $k(\alpha)$. The total ordering of
$B(\alpha)$ is denoted by $\leq_\alpha$. Sometimes we write
$\alpha=(\mathfrak{p}_\alpha,\leq_\alpha)$.
\begin{deft}\label{sper} The real spectrum of $B$, denoted by Sper $B$, is the
  collection of all pairs $\alpha=(\mathfrak{p}_\alpha,\leq_\alpha)$,
  where $\mathfrak{p}_\alpha$ is a prime ideal of $B$ and
  $\leq_\alpha$ is a total ordering of $B/\mathfrak{p}_{\alpha}$.
\end{deft}
Given a point $\delta \in \sper(A)$ and an element $f\in A$, the notation $|f|_\delta$ will mean $f$ if $f\ge_\delta0$,
$-f$ if $f\le_\delta0$. When no confusion is possible, we will write simply $|f|$, with $\delta$ understood.\medskip

Two kinds of points occur in $\sper\ B$: finite points and points at infinity.
\begin{deft}\label{bounded} Let $B$ be an $R$-algebra and $\alpha$ a
  point of $\sper\ B$. We say that $\alpha$ is {\bf finite} if for
  each $y\in B[\alpha]$ there exists $N\in R$ such that
  $|y|_\alpha<_\alpha N$. Otherwise, we say that $\alpha$ is {\bf a point at infinity}.
\end{deft}

\noi{\bf Notation:} The subset of $\sper\ B$ consisting of all the finite
points will be denoted by $\sper^*B$.\medskip

It is known, as we explain in detail in \S\ref{Cu}, that $\sper\ B$ is closely related to the
space $\bigcup\limits_{\mathfrak p\in\mbox{Spec}\ B}S_{\mathfrak p}$, where $S_{\mathfrak p}$ denotes the
Zariski--Riemann surface of the residue field $\kappa(\mathfrak p)$. Namely, one can associate to every point
$\delta\in\sper\ B$ a valuation $\nu_\delta$ of $\kappa(\mathfrak p_\delta)$ (where $\mathfrak p_\delta$ is
the support of $\delta$) with totally ordered residue field $k_\delta$. Conversely, given a prime ideal
$\mathfrak{p}\subset B$ and a valuation $\nu$ of $\kappa(\mathfrak p)$ with totally ordered residue field,
one can define a point $\delta\in\sper\ R[x_1,\dots,x_n]$ with $\mathfrak p_\delta=\mathfrak p$ and $\nu_\delta=\nu$ by
specifying the signs of finitely many elements of $\kappa(\mathfrak p)$ with respect to the total ordering $\le_\delta$ (see
Remark \ref{BaerKrull} below).

The real spectrum $\sper\ B$ is endowed with the {\bf spectral (or Harrison) topology}. By
definition, this topology has basic open sets of the form
$$
U(f_1,\ldots,f_k)=\{\alpha\ |\ f_1(\alpha)>0,\ldots,f_k(\alpha)>0\}
$$
with $f_1,...,f_k\in B$. Here and below, we commit the following
standard abuse of notation: for an element $f\in B$, $f(\alpha)$
stands for the natural image of $f$ in $B[\alpha]$ and the inequality
$f(\alpha)>0$ really means $f(\alpha)>_\alpha0$.

Denote by Maxr$(A)$ the set of points $\alpha \in \sper(A)$ such that
$\mathfrak{p}_\alpha$ is a maximal ideal of $A$. We view Maxr$(A)$ as
a topological subspace of $\sper(A)$ with the spectral (respectively,
constructible) topology. We may naturally identify $R^n$
with Maxr$(A)$: a point $(a_1,\ldots,a_n)\in R^n$ corresponds to the point
$\alpha=(\mathfrak{p}_\alpha,\leq_\alpha)\in\sper(A)$, where $\mathfrak{p}_\alpha$
is the maximal ideal
$$
\mathfrak{p}_\alpha=(x_1-a_1,\ldots,x_n-a_n)
$$
and $\leq_\alpha$ is the unique order on $R$. The spectral topology on $\sper(A)$ induces
the euclidean topology on $R^n$.

Let $A=R[x_1,\dots,x_n]$. Take a point $\delta\in\sper\ A$. In \S\ref{infinity} we
associate to $\delta$ three disjoint subsets $I_\delta,F_\delta,G_\delta\subset\{1,\dots,n\}$, as follows. By definition,
the set $I_\delta\coprod F_\delta\coprod G_\delta$ is the set of all $j\in\{1,\dots,n\}$ such that
\begin{equation}
\nu_\delta(x_j)=0.\label{eq:nudelta0}
\end{equation}
The set $G_\delta$ consists of all $j$ such that $|x_j|_\delta$ is bounded below by all the elements
of $R$, $I_\delta$---of all $j$ such that (\ref{eq:nudelta0}) holds and $|x_j|_\delta$ is smaller than any strictly positive
constant in $R$ and $F_\delta$ of all $j$ such that $|x_j|_\delta$ is bounded both above and below by strictly positive
constants from $R$. We show that $I_\delta=\emptyset$ whenever $G_\delta=\emptyset$.

Let $T$ be a set such that $G_\delta\subset T\subset G_\delta\cup F_\delta$. For $j\in\{1,\dots,n\}$, let $y_j=x_j$ if
$j\notin T$ and $y_j=\frac1{x_j}$ if $j\in T$. Let
$B_T=R[y_1,\dots,y_n]$. We associate to $\delta$ a point $\delta^*$ in $\sper^*B_T$ such that
$A(\delta)=B_T(\delta^*)$. We show that $R_\delta$ is a localization of $R_{\delta^*}$ at a prime ideal.

Let $I,F,G$ be three disjoint subsets of $\{1,\dots,n\}$, such that if $G=\emptyset$ then $I=\emptyset$. Let $U_{I,F,G}$
denote the set of all points of $\sper\ A$ such that $I=I_\delta$, $F=F_\delta$ and $G=G_\delta$. The main theorem, Theorem
\ref{homeo}, describes a homeomorphism between $U_{I,F,G}$ and a certain explicitly described subspace
$U^*_{I,F,G}\subset\sper^*(B_T)$, where $T$ is a set satisfying $G\subset T\subset F\cup G$. At the end of section
\S\ref{infinity}, we describe a partition
\begin{equation}
\sper(A)=\coprod\limits_{I,F,G}U_{I,F,G},\label{eq:partition}
\end{equation}
where $I,F,G$ runs over all the triples of disjoint subsets of $\{1,\dots,n\}$
such that $I=\emptyset$ whenever $G=\emptyset$, and each $U_{I,F,G}$ is homeomorphic to a
subspace $U^*_{I,F,G}\subset\sper^*(B_T)$, as above.

This paper originally grew out of the authors' joint work with J.J. Madden \cite{LMSS} on the Pierce--Birkhoff conjecture.
Certain definitions and constructions only worked for finite points of $\sper\ A$, so a need naturally arose to cover
$\sper A$ by subspaces, each of which is homeomorphic to a subspace of $\sper^*B$ for some other polynomial ring $B$.
Eventually, we found another way of getting around this difficulty and were able to deal in a uniform way with all the
points of $\sper\ A$, whether finite or infinite. However, we hope that the decomposition (\ref{eq:partition}) may
some day come in useful to someone who is faced with finiteness problems similar to ours. Also, since in \cite{LMSS} we
are interested in proving connectedness of certain subsets of $\sper\ A$, we gave a variation of the decomposition
(\ref{eq:partition}) into sets which are \textit{not} disjoint; we derive it as an easy consequence
of (\ref{eq:partition}).

\section{The valuation associated to a point in the real spectrum}
\label{Cu}

Let $B$ be a ring and $\alpha$ a point in $\sper\ B$. In this section
we define the valuation $\nu_\alpha$ of $B(\alpha)$, associated to
$\alpha$. We also give a geometric interpretation of points in the real spectrum as semi-curvettes.
\medskip

\noi{\bf Terminology}: If $B$ is an integral domain, the phrase ``valuation
of $B$'' will mean ``a valuation of the field of fractions of $B$, non-negative on $B$''. Also, we will sometimes commit
the following abuse of notation. Given a ring $B$, a prime ideal $\mathfrak p\subset B$, a
valuation $\nu$ of $\frac B{\mathfrak p}$
and an element $x\in B$, we will write $\nu(x)$ instead of $\nu(x\mod\mathfrak p)$, with the usual convention that
$\nu(0)=\infty$, which is taken to be greater than any element of the value group.
\medskip

For a point $\alpha$ in $\sper\ B$, we define the valuation ring $R_\alpha$ by
$$
R_\alpha=\{x\in B(\alpha)\ |\ \exists z\in B[\alpha],|x|_\alpha\leq_\alpha z\}.
$$
That $R_\alpha$ is, in fact, a valuation ring, follows because for any
$x\in B(\alpha)$, either $x\in R_\alpha$ or $\frac1x\in R_\alpha$. The
maximal ideal of $R_\alpha$ is $M_\alpha=\left\{x\ \in B(\alpha)
  \left|\ |x|_\alpha<\frac{1}{|z|_\alpha}, \ \forall z\in B[\alpha] \setminus \{0\}
  \right.\right\}$; its residue field $k_\alpha$ comes equipped
with a total ordering, induced by $\le_\alpha$. For a ring $B$ let $U(B)$ denote the multiplicative group of units of $B$.
Recall that $\dsp{\Gamma_\alpha \cong\frac{B(\alpha)\setminus \{0\}}{U(R_\alpha)}}$ and that the
valuation $\nu_\alpha$ can be identified with the natural homomorphism \\
\centerline{$\dsp{B(\alpha)\setminus \{0\} \to \frac{B(\alpha)\setminus
  \{0\}}{U(R_\alpha)}}$.}
\medskip

By definition, we have a natural ring homomorphism
\begin{equation}
B\rightarrow R_\alpha\label{eq:hom}
\end{equation}
whose kernel is $\mathfrak{p}_\alpha$. The valuation $\nu_\alpha$ has the following properties:

(1) $\nu_\alpha(B[\alpha]) \geq 0$

(2) If $B$ is an $R$-algebra then for any positive elements $y,z \in B(\alpha)$,
\begin{equation}
\nu_\alpha(y) < \nu_\alpha(z) \implies y > Nz,\ \forall N \in
R \label{eq:val2}
\end{equation}
(an example at the end of the paper shows that the converse implication in (\ref{eq:val2}) is not true in
general).

\begin{rek} Let $B$ be an $R$-algebra and take a point $\alpha\in\sper^*B$ (see Definition
\ref{bounded}). Then
\begin{equation}
R_\alpha=\left\{x\in B(\alpha)\ \left|\ \exists N\in R,|x|\leq_\alpha
N\right.\right\}.\label{eq:myvaluation}
\end{equation}
Thus for points $\alpha\in\sper^*B$ the valuation $\nu_\alpha$ of $B(\alpha)$
depends on the ordering $\le_\alpha$ but not on the ring $B[\alpha]$ (this means that given another $R$-algebra $\tilde B$,
a point $\tilde\alpha\in\sper^*\tilde B$ and an order-preserving isomorphism $\phi:B(\alpha)\cong\tilde B(\tilde\alpha)$, we
have $\phi(R_\alpha)=R_{\tilde\alpha}$).
\end{rek}
\begin{rek}\label{BaerKrull} (\cite{Baer}, \cite{Krull}, \cite{BCR} 10.1.10,
p. 217) Conversely, the point $\alpha$ can be reconstructed
from the ring $R_\alpha$ by specifying a certain number of sign
conditions (finitely many conditions when $B$ is noetherian), as we
now explain. Take a prime ideal $\mathfrak{p}\subset B$ and a
valuation $\nu$ of $\kappa({\mathfrak{p}}):=\frac{B_{\mathfrak{p}}}{{\mathfrak{p}}B_{\mathfrak{p}}}$,
with value group $\Gamma$. Let
$$
r=\dim_{\F_2}(\Gamma/2\Gamma)
$$
(if $B$ is not noetherian, it may happen that $r=\infty$). Let
$x_1,\ldots,x_r$ be elements of $\kappa({\mathfrak{p}})$ such that
$\nu(x_1),\ldots,\nu(x_r)$ induce a basis of the $\F_2$-vector space
$\Gamma/2\Gamma$. Then for every $x\in\kappa({\mathfrak{p}})$, there exist $f\in\kappa({\mathfrak{p}})$ and a unit $u$ of
$R_\nu$ such that $x=ux_1^{\epsilon_1}\cdots x_r^{\epsilon_r}f^2$ with $\epsilon_i\in\{0,1\}$
(to see this, note that for a suitable choice of $f$ and $\epsilon_j$ the
value of the quotient $u$ of $x$ by the product $x_1^{\epsilon_1}\cdots
x_r^{\epsilon_r}f^2$ is 0, hence $u$ is invertible in $R_\nu$). Now,
specifying a point $\alpha\in\sper\ B$ supported at ${\mathfrak{p}}$ amounts
to specifying a valuation $\nu$ of $\frac B{\mathfrak{p}}$, a total ordering of the residue
field $k_\nu$ of $R_\nu$, and the sign data $\sgn\ x_1,\dots,\sgn\ x_r$. For $x\notin\mathfrak{p}$, the sign of $x$ is
given by the product $\sgn(x_1)^{\epsilon_1}\cdots\sgn(x_r)^{\epsilon_r}\sgn(u)$, where $\sgn(u)$
is determined by the ordering of $k_\nu$.
\end{rek}

Points of $\sper\ B$ admit the following geometric interpretation (we refer the reader to \cite{Fuc}, \cite{Kap}, \cite{Pre},
p. 89 and \cite{PC} for the construction and properties of generalized power series rings and fields).

\begin{deft} Let $k$ be a field and $\Gamma$ an ordered abelian group. The
  generalized formal power series ring $k\left[\left[t^\Gamma\right]\right]$
  is the ring formed by elements of the form $\sum\limits_\gamma a_\gamma
  t^\gamma$, $a_\gamma\in k$ such that the set $\left\{\gamma\ \left|\ a_\gamma\ne0\right.\right\}$
  is well ordered.
\end{deft}
The ring $k\left[\left[t^\Gamma\right]\right]$ is equipped with the
natural $t$-adic valuation $v$ with values in $\Gamma$, defined by
$v(f)=\inf\{\gamma\ |\ a_\gamma\neq0\}$ for $f=\sum\limits_\gamma
a_\gamma t^\gamma\in k\left[\left[t^\Gamma\right]\right]$. Specifying
a total ordering on $k$ and $\dim_{\F_2}(\Gamma/2\Gamma)$ sign
conditions defines a total ordering on
$k\left[\left[t^\Gamma\right]\right]$. In this ordering $|t|$ is
smaller than any positive element of $k$. For example, if $t^\gamma>0$ for all
$\gamma\in\Gamma$ then $f>0$ if and only if $a_{v(f)}>0$.

For an ordered field $k$, let $\bar k$ denote the real closure of $k$.
The following result is a variation on a theorem of Kaplansky (\cite{Kap}, \cite{Kap2})
for valued fields equipped with a total ordering.

\begin{thm} \textnormal{\textbf{(\cite{PC}, p. 62, Satz 21)}} Let $K$ be a
  real valued field, with residue field $k$ and value group $\Gamma$.
  There exists an injection $K\hookrightarrow\bar
  k\left(\left(t^\Gamma\right)\right)$ of real valued fields.
\end{thm}
\medskip

Let $\alpha\in\sper\ B$ and let $\Gamma_\alpha$ be the value group of
$\nu_\alpha$. In view of (\ref{eq:hom}) and the Remark above, specifying a
point $\alpha\in\sper\ B$ is equivalent to specifying a total order of $k_\alpha$, a morphism
$$
B[\alpha]\to\bar k_\alpha\left[\left[t^{\Gamma_\alpha}\right]\right]
$$
and $\dim_{\F_2}(\Gamma_\alpha/2\Gamma_\alpha)$ sign conditions as above.\medskip

\noi We may pass to usual spectra to obtain morphisms
$$
\mbox{Spec}\ \left(\bar k_\alpha\left[\left[t^{\Gamma_\alpha}\right]\right]\right)\to\mbox{Spec}\
B[\alpha]\to\mbox{Spec}\ B.
$$
In particular, if $\Gamma_\alpha=\Z$, we obtain a
\textbf{formal curve} in Spec $B$ (an analytic curve if the series are
convergent). This motivates the following definition:

\begin{deft} Let $k$ be an ordered field. A $k$\textbf{-curvette} on $\sper(B)$ is a morphism of the form
$$
\alpha:B\to k\left[\left[t^\Gamma\right]\right],
$$
where $\Gamma$ is an ordered group. A $k$\textbf{-semi-curvette} is a $k$-curvette $\alpha$ together with
  a choice of the sign data $\sgn\ x_1$,..., $\sgn\ x_r$, where
  $x_1,...,x_r$ are elements of $B$ whose $t$-adic values induce an
  $\F_2$-basis of $\Gamma/2\Gamma$.
\end{deft}

We have thus explained how to associate to a point $\alpha$ of Sper\ $B$ a
$\bar k_\alpha$-semi-curvette. Conversely, given an
ordered field $k$, a $k$-semi-curvette $\alpha$ determines a prime
ideal $\mathfrak{p}_\alpha$ (the ideal of all the elements of $B$
which vanish identically on $\alpha$) and a total ordering on
$B/\mathfrak{p}_\alpha$ induced by the ordering of the ring
$k\left[\left[t^\Gamma\right]\right]$ of formal power series. These
two operations are inverse to each other. This establishes a
one-to-one correspondence between semi-curvettes and points of $\sper\
B$.\medskip

Below, we will sometimes describe points in the real spectrum by specifying the corresponding semi-curvettes.\medskip

\noi{\bf Example:} Consider the curvette $R[x,y]\to R[[t]]$
defined by $x\mapsto t^2,\ y\mapsto t^3$, and the semi-curvette given
by declaring, in addition, that $t$ is positive. This semi-curvette is nothing but the upper
branch of the cusp.\medskip

Later in the paper, we will need, for a certain number $p\in\{0,1,\dots,n\}$ and two points $\delta,\delta^*$ living in
different spaces, to compare $(n-p)$-tuples of elements
such as $(\nu_\delta(x_{p+1}),\dots,\nu_\delta(x_n))\in\Gamma_\delta^{n-p}$ and
$(\nu_{\delta^*}(y_{p+1}),\dots,\nu_{\delta^*}(y_n))\in\Gamma_{\delta^*}^{n-p}$ and to be able to say that they are in
some sense ``equivalent''. To do this, we need to embed $\Gamma_\delta$ in some ``universal'' ordered group.\medskip

\noi{\bf Notation and convention:} Let us denote by $\Gamma$ the ordered group $\R^n_{lex}$. This means that elements of
$\Gamma$ are compared as words in a dictionary: we say that $(a_1,\dots,a_n)<(a'_1,\dots,a'_n)$ if and only if there
exists $j\in\{1,\dots,n\}$ such that $a_q=a'_q$ for all $q<j$ and $a_j<a'_j$.

The reason for introducing $\Gamma$ is that by
Abhyankar's inequality we have $\mbox{rank}\ \nu_\delta\le\dim A=n$ for all $\delta\in\sper\ A$, so the value group
$\Gamma_\delta$ of $\nu_\delta$ can be embedded into $\Gamma$ as an ordered subgroup (of course, this embedding is far from
being unique). Let $\Gamma_+$ be the semigroup of non-negative elements of $\Gamma$. \medskip

Fix a strictly positive integer $\ell$. In order to deal rigourously
with $\ell$-tuples of elements of $\Gamma_\delta$ despite the non-uniqueness of the embedding $\Gamma_\delta\subset\Gamma$,
we introduce the category $\mathcal{O}\mathcal{G}\mathcal{M}(\ell)$, as follows. An
object in $\mathcal{O}\mathcal{G}\mathcal{M}(\ell)$ is an ordered
abelian group $G$ together with $\ell$ fixed
generators $a_1,\dots,a_\ell$ (such an object will be denoted by
$(G,a_1,\dots,a_\ell)$). A morphism from $(G,a_1,\dots,a_\ell)$
to $(G',a'_1,\dots,a'_\ell)$ is a homomorphism $G\rightarrow G'$ of
ordered groups which maps $a_j$ to $a'_j$ for each $j$.
\medskip

Given $(G,a_1,\ldots,a_\ell),\ (G',a'_1,\ldots,a'_\ell) \in
Ob(\mathcal{OGM}(\ell))$, the notation
\begin{equation}
(a_1,\ldots,a_\ell) \underset{\circ}{\sim} (a'_1,\ldots,a'_\ell)
\end{equation} will mean that $(G,a_1,\ldots,a_\ell)$ and
$(G',a'_1,\ldots,a'_\ell)$ are isomorphic in $\mathcal{OGM}(\ell)$.
\medskip

Take an element
$$
a=(a_1,\ldots,a_\ell)\in\Gamma_+^\ell.
$$
Let $G\subset\Gamma$ be the ordered group generated by
$a_1,\ldots,a_\ell$. Then $(G,a_1,\ldots,a_\ell) \in
Ob(\mathcal{OGM}(\ell))$. For each $\delta \in \sper(A)$, let
$\Gamma_\delta$ denote the value group of the associated valuation
$\nu_\delta$ and $\Gamma_\delta^*$ the subgroup of $\Gamma_\delta$
generated by $\nu_\delta(x_1)$, \dots, $\nu_\delta(x_n)$. In this way,
we associate to $\delta$ an object $\left(\Gamma_\delta^*,
  \nu_\delta(x_1), \dots, \nu_\delta(x_n)\right)\in
Ob(\mathcal{OGM}(n))$.
\medskip

\noi{\bf Notation.} Let $\Gamma$ be an ordered group. Consider an
$\ell$-tuple $a=(a_1,\ldots,a_\ell)\in\Gamma^\ell$. We denote by
$Rel(a)$ the set
$$
Rel(a)=\left\{(m_1,\dots,m_\ell,m_{\ell+1},\dots,m_{2\ell})\in
  \Z^{2\ell}\ \left|\ \sum\limits_{j=1}^\ell m_ja_j>0 \mbox{ and }
    \sum\limits_{j=\ell+1}^{2\ell} m_ja_{j-\ell}=0\right.\right\}.
$$
\medskip

\begin{rek}\label{Rel=ogm} Let $\Gamma$ and $a$ be as above and let $G$ be the
subgroup of $\Gamma$ generated by $a_1,\ldots,a_\ell$, so that
$(G,a_1,\ldots,a_\ell) \in Ob(\mathcal{OGM}(\ell))$. The
set $Rel(a)$ completely determines the isomorphism class of
$(G,a_1,\ldots,a_\ell)$ in $\mathcal{OGM}(\ell)$ and vice-versa; the
set $Rel(a)$ and the isomorphism class of
$(G,a_1,\ldots,a_\ell)$ are equivalent sets of data.
\end{rek}

\section{Points at infinity of \sper(A)}
\label{infinity}

 In this section, we study the structure of the set of points at infinity
in $\sper(A)$. Take a point $\delta \in \sper(A)$. Renumbering the coordinates if necessary, we may assume
 there exists $p$, $0\leq p \leq n$, such that
 \begin{equation}
 \nu_\delta(x_j)=0\text{ for }
 1 \leq i \leq p\text{ and }\nu_\delta(x_j)>0\text{ for }j > p.\label{eq:defp}
 \end{equation}
 For a subset $T$ of $\{1,\ldots,p\}$, let $B_T=R[y_1,\dots,y_n]$ where
 \begin{eqnarray}
 y_j&=&x_j\qquad\text{if }j\in\{1,\dots,n\}\setminus T \label{yj}\\
  &=& 1/x_j \quad\text{if }j\in T \label{y1j}.
 \end{eqnarray}
For certain subsets $T\subset\{1,\ldots,p\}$ we will associate to $\delta$ a point $\delta^*$ in $\sper^*(B_T)$
 such that $A(\delta)=B_T(\delta^*)$. We will define a new valuation
 $\nu_{\delta^*}$ of $A(\delta)$, such that $R_\delta$ is a
 localization of $R_{\delta^*}$ at a suitable prime ideal.
 At the end of this section, we will use these results to cover
 $\sper(A)$ by sets, each of which is homeomorphic to a certain
 subspace of $\sper^*(B_{T})$ for some $T$.
 \medskip

 First, take any subset $T \subset \{1,\ldots,n\}$ whatsoever. Let
 $B_T$ be defined as in (\ref{yj})--(\ref{y1j}).
 \medskip

 \noi \textbf{Notation }: The notation $A_f$ stands for the
 localization of $A$ by $f$, the ring $A[1/f]$.
 \medskip

 \noi \textbf{Remark} : We have a natural homeomorphism
 \begin{equation}
 \sper(A) \setminus \{f=0\} \stackrel{\sim}{\rightarrow} \sper(A_f)
 \end{equation}
 \medskip

 \noi Consider the natural isomorphism $A_{\prod{x_j},j \in T} \cong
 (B_T)_{\prod{y_j}, j \in T}$. It induces a homeomorphism
 \begin{equation}
 \begin{array}{cccc}
 \psi:& \sper(A) \setminus \left\{\prod_{j \in T}x_j=0\right\}
   & \to & \sper(B_T) \setminus \left\{\prod_{j \in T}y_j=0\right\} \\ & \parallel
   & & \parallel \\ & \sper(A_{\prod_{j \in T}x_j}) & \to &
   \sper((B_T)_{\prod_{j \in T}y_j})
 \end{array}
 \end{equation}
 which we describe explicitly for future reference.
Take a point $\delta \in \sper(A_{\prod_{j \in T}x_j})$. We will
 now describe the point $\delta^* = \psi(\delta)$ in $\sper(B_T)$, as
 follows. The ideal
 $\mathfrak{p}_{\delta^*}$ is the prime ideal of $B_T$ such
 that
 \begin{equation} \mathfrak{p}_\delta A_{\prod{x_j},j \in T}
   \cong \mathfrak{p}_{\delta^*}(B_T)_{\prod{y_j}, j \in
 T} \label{decadix}
 \end{equation}
 Then (\ref{decadix}) implies the existence of a canonical isomorphism
 \begin{equation}
  \phi: \kappa(\mathfrak{p}_\delta) \cong \kappa(\mathfrak{p}_{\delta^*})
   \label{decadix2}.
   \end{equation}
   The total order $\leq_{\delta^*}$ is the order induced by $\delta$ on $\kappa(\mathfrak{p}_{\delta^*})$ via the
 isomorphism (\ref{decadix2}). This describes $\psi$; the inverse map $\psi^{-1}$ is described in a completely
 analogous way.

Here and below, $R_{>0}$ will denote the set of strictly positive elements of $R$.

Take a $\delta \in \sper(A)$ and let $p$ be as in (\ref{eq:defp}). We associate to $\delta$ a partition
 $$
 \{1,\ldots,p\} = I_\delta \coprod F_\delta \coprod G_\delta,
 $$
 as follows:
\begin{eqnarray}
j &\in& I_\delta \Longleftrightarrow |x_j|_\delta <_\delta \epsilon, \ \forall \epsilon \in
 R_{>0}\label{eq:Idelta}\\
 j& \in& F_\delta \Longleftrightarrow \exists c_1,c_2 \in R_{>0}\quad\text{such that
}c_1<_\delta |x_j|_\delta<_\delta c_2\label{eq:Fdelta}\\
j &\in& G_\delta \Longleftrightarrow |x_j|_\delta >_\delta N, \ \forall N \in R.\label{eq:Gdelta}
\end{eqnarray}
 \begin{rek} We have $\delta \in \sper^*(A)$ if and only if
 $G_\delta= \emptyset$. Below, we show that in this case necessarily $I_\delta=\emptyset$.
 \end{rek}
 \noi Take a set $T$ such that
\begin{equation}
 G_\delta\subset T\subset G_\delta\cup F_\delta.\label{eq:Tdef}
\end{equation}
Let $B_T$ be the ring defined by (\ref{yj}) and  (\ref{y1j}). It follows from (\ref{eq:Fdelta}), (\ref{eq:Gdelta}) and
(\ref{eq:Tdef}) that $x_j \notin \mathfrak{p}_\delta$ for $j \in T$. Let $\delta^*=\psi(\delta)$. It is immediate from
 the definition  that $\delta^*$ is finite in $\sper(B_T)$.

 \begin{prop}\label{propdelta} The valuation $\nu_{\delta^*}$ of
   $B_T(\delta^*)$ associated to $\delta^*$ has the following
   properties:

 (1) $\nu_{\delta^*}(y_j) = 0$ for $j \in F_\delta$;

 (2) $\nu_{\delta^*}(y_j) > 0$ for $j \in I_\delta \cup G_\delta$;

 (3) there exists $q \in G_\delta$ and a strictly positive integer $N$
 such that, for all $j \in I_\delta$,
 \begin{equation}
 N\nu_{\delta^*}(y_q)> \nu_{\delta*}(y_j).\label{eq:Nq>j}
 \end{equation}
 In particular, if $I_\delta\ne\emptyset$ then $G_\delta\ne\emptyset$.

 (4) The valuation ring $R_\delta$ is the localization of $R_{\delta^*}$
 at a prime ideal; this gives rise to a surjective order-preserving group homomorphism
 $\tilde{\phi}: \Gamma_{\delta^*} \to \Gamma_\delta$ whose kernel is an
 isolated subgroup.

 (5) For all $j\in\{1,\dots,n\}$, $\tilde{\phi}(\nu_{\delta^*}(y_j)) =
 \nu_\delta(x_j)$.

 (6) For $j \in \{1,\ldots,p\}$, $\nu_{\delta^*}(y_j)
 \in \ker(\tilde{\phi})$. In particular, given any $j \in \{1,\ldots,p\}, \ t
 \in \{p+1,\ldots,n\}$ and $N' \in \N$, we have $N' \nu_{\delta^*}(y_j) <
 \nu_{\delta^*}(y_t)$.

 (7) Assume that $\nu_\delta(x_{p+1}),\ldots,\nu_\delta(x_n)$ are $\Q$-linearly independent.
 Then
 $$
 (\nu_{\delta^*}(y_{p+1}),\ldots,\nu_{\delta^*}(y_n)) \underset{\circ}{\sim}
  (\nu_{\delta}(x_{p+1}),\ldots,\nu_{\delta}(x_n))
 $$
 in $\mathcal{OGM}(n-p)$.
 \end{prop}

 \noi Proof : (1) Take $j\in F_\delta$. We have $1/|y_j|_{\delta^*} <_{\delta^*} c$ for some
 $c\in R$ by definition of $F_\delta$ (\ref{eq:Fdelta}). Hence $\frac1{y_j}\in R_{\delta^*}$ and the result follows.

 (2) Take $j \in I_\delta \cup G_\delta$. By definition of $I_\delta$ (\ref{eq:Idelta}),
$G_\delta$ (\ref{eq:Gdelta}) and (\ref{eq:Tdef}), $|y_j|_{\delta^*} <_{\delta^*} \epsilon$ for every $\epsilon\in
 R_{>0}$, so $1/|y_j|_{\delta^*} >_{\delta^*} N$ for every $N \in R$. By the boundedness of
 $\delta^*$, for each $f \in B_T$, we have
 $|f(y_1,\ldots,y_n)|_{\delta^*} <_{\delta^*} N'$ for some $N' \in R$. Hence $1/|y_j|_{\delta^*} >_{\delta^*}
 f(y_1,\ldots,y_n)$ for each $f \in B_T$, so $1/y_j \notin
 R_{\delta^*}$. This proves that $\nu_{\delta^*}(y_j) > 0$.

 (3) Take a $j \in I_\delta$. Since $\nu_\delta(x_j)=0$, we have $1/x_j \in R_\delta$. This
 means that there exists $z \in A[\delta]$ such that $1/|x_j|_\delta <_\delta |z|_\delta$. Now,
 $z$ is a polynomial in the $x_k$, $k=1,\ldots,n$, and taking
 $x_q$, $q \in G_\delta$, such that $|x_q|_\delta \geq_\delta |x_k|_\delta$ for all $k
 \in G_\delta$ (and hence $|x_q|_\delta \geq_\delta |x_k|_\delta$ for all $k \in
 \{1,\ldots,n\}$), there exists $N >0$ such that $1/|x_j|_\delta <_\delta |x_q|_\delta^N$
 for all $j \in I_\delta$, so that $|y_j|_{\delta^*} >_{\delta^*} |y_q|_{\delta^*}^N$. Then
 $N\nu_{\delta^*}(y_q) \geq \nu_{\delta^*}(y_j)$ by equation (\ref{eq:val2}). Replacing $N$ by $N+1$, we can make the
 inequality (\ref{eq:Nq>j}) strict.

 (4) It is well known that every homomorphism between two valuation
 rings having the same field of fractions is a localization at a prime
 ideal. Thus it is sufficient to show that $R_{\delta^*} \subset
 R_\delta$. Take $f \in R_{\delta^*}$. By definition, this means that
 $|f|_{\delta^*}$ is bounded above by a polynomial in the $y_j$ with respect to
 $\leq_{\delta^*}$, and hence also by a monomial $\omega$ in the $y_j$. Then
 $\phi^{-1}(\omega)$ is bounded above with
 respect to $\leq_\delta$ by a monomial in the $x_j$, in which $x_j$ with $j\in T$ appear with non-positive exponents.
Since each $\frac1{|x_j|_\delta}$, $j\in T$, is bounded above by a constant in $R$, replacing factors of the form
$x_j^{-\gamma_j}$, $j\in T$, $\gamma_j\in\mathbb N$ by a suitable constant in $R$, we obtain that
$\phi^{-1}(f)$ is bounded above with respect to $\le_{\delta^*}$ by a monomial in $y$ with non-negative exponents. This
proves that $\phi^{-1}(f) \in R_\delta$ as desired.

The last statement of (4) follows immediately by
 the general theory of composition of valuations (\cite{ZS}, Chapter VI, \S10, p. 43). Alternatively, recall that
 $\dsp{\Gamma_\delta \cong\frac{A(\delta)\setminus \{0\}}{U(R_\delta)}}$ and that the
 valuation $\nu_\delta$ can be identified with the natural homomorphism \\
 \centerline{$\dsp{A(\delta)\setminus \{0\} \to \frac{A(\delta)\setminus
   \{0\}}{U(R_\delta)}}$.} Similarly, $\nu_{\delta^*}$ can be thought
 of as \\ \centerline{$\dsp{B_T(\delta^*) \setminus \{0\} \to
 \frac{B_T(\delta^*) \setminus \{0\}}{U(R_{\delta^*})} \cong
   \Gamma_{\delta^*}}$.} From the isomorphism $\phi$ and the
   inclusion $R_{\delta^*} \hookrightarrow R_\delta$, we obtain a
   natural surjective homomorphism of ordered groups
   \begin{equation}
     \tilde{\phi} : \frac{B_T(\delta^*) \setminus
       \{0\}}{U(R_{\delta^*})} \to \frac{A(\delta) \setminus
       \{0\}}{U(R_\delta)}.
       \end{equation}

 (5) If $j \notin T$, the fact that $\phi(x_j)= y_j$ implies that
$$
\tilde{\phi}(y_j\mod U(R_{\delta^*})) = x_j \mod U(R_\delta).
$$
If $j \in T$, we have $\phi(x_j)= 1/y_j$, hence
 $$
 \tilde{\phi}(\nu_{\delta^*}(y_j)) = \nu_\delta(1/x_j) = 0 =\nu_\delta(x_j).
 $$

 (6) is an immediate consequence of (5) and the fact that
 $\nu_\delta(x_1)=\cdots=\nu_\delta(x_p)=0$.

 (7) By Remark \ref{Rel=ogm} at the end of the previous section, it suffices to
 prove that
 \begin{equation}
 \label{de_rel}
   Rel(\nu_{\delta^*}(y_{\ell+1}),\ldots,\nu_{\delta^*}(y_n))
 =
 Rel(\nu_\delta(x_{\ell+1}),\ldots,\nu_{\delta}(x_n)).
 \end{equation}

 The fact that $\nu_\delta(x_{p+1}),\ldots,\nu_\delta(x_n)$ are $\Q$-linearly
 independent and (5) of the Proposition imply that so are
 $\nu_{\delta^*}(y_{p+1})$, \ldots, $\nu_{\delta^*}(y_n)$. Hence, using (5) of the Proposition
 again, for any $(n-p)$-tuple, $(m_{p+1},\ldots,m_n) \in
   \Z^{n-p}$, we have $\sum_{i=p+1}^n m_j\nu_\delta(x_j) > 0$ if and only if
   $\sum_{j=p+1}^n m_j\nu_{\delta^*}(y_j) > 0$. Together with the
   linear independence of $\nu_\delta(x_{p+1}),\ldots,\nu_\delta(x_n)$ and of
   $\nu_{\delta^*}(y_{p+1}),\ldots,\nu_{\delta^*}(y_n)$, this proves
   the desired equality (\ref{de_rel}). \hfill $\Box$
 \medskip

 Let $G$ be an ordered group of rank $r$ and $\ell$ a positive integer. Take $\ell$ elements $a_1,\ldots,a_\ell \in G$. Let
 $(0) = \Delta_r \subsetneqq \Delta_{r-1} \subsetneqq \cdots \subsetneqq
 \Delta_0 =G$ be the isolated subgroups of $G$. Renumbering the $a_j$ if
 necessary, we may assume that there exist integers $i_0,i_1,\ldots,i_r$ with
 $$
 \ell =i_0 \geq i_1 \geq i_2 \geq \cdots\geq i_r = 0,
 $$
 such that $a_{i_{q+1}}, \ldots, a_{i_q} \in \Delta_q-\Delta_{q+1}$ for $q \in
 \{0,\ldots,r-1\}$.

 \begin{deft}\label{scalewise}
 We say that $a_1,\ldots,a_\ell$ are scalewise $\Q$-linearly
 independent if, for each $q \in \{0,\ldots,r-1\}$, the images of
 $a_{i_{q+1}},\ldots,a_{i_q}$ in $\dsp{\frac{\Delta_q}{\Delta_{q+1}}}$
 are $\Q$-linearly independent.
 \end{deft}
 \begin{rek} Let the notation be as above and assume that
 $a_1,\ldots,a_\ell$ are scalewise $\Q$-linearly independent. Let
 $\lambda:G\rightarrow G'$ be a homomorphism of ordered groups. Then
 $\lambda(a_1),\ldots,\lambda(a_\ell)$ are scalewise $\Q$-linearly
 independent if and only if they are $\Q$-linearly independent if and
 only if all of them are non-zero. This is precisely the form in which
 we will use scalewise $\Q$-linear independence in the sequel.
 \end{rek}
 Fix an integer $p \in \{1,\ldots,n\}$ and two decompositions
\begin{equation}
\{1,\ldots,p\} = H \coprod T=I\coprod F\coprod G,\label{eq:decomp}
\end{equation}
where $I=\emptyset$ whenever $G=\emptyset$,
\begin{eqnarray}
I&\subset&H\quad\text{and}\label{eq:IdansH}\\
G&\subset&T.\label{eq:GdansT}
\end{eqnarray}
Fix an $n$-tuple $(a_1,\ldots,a_n) \in\Gamma_+^n$ such that $a_1=\cdots = a_p=0$. Let
\begin{equation}
 U_{I,F,G} = \left\lbrace \delta \in \sper(A) \ \left|
 \begin{array}{l}
\forall j \in I, \forall c \in R_{>0},\ |x_j|_\delta <_\delta c \\
\forall j\in F \exists c_1,c_2\in R_{>0},c_1<_\delta|x_j|_\delta<_\delta c_2\\
\forall j \in G, \forall N \in R,\ |x_j|_\delta >_\delta N  \\
 \nu_\delta(x_1)=\ldots=\nu_\delta(x_p)=0\\
\nu_\delta(x_{p+1})>0,\ldots,\nu_\delta(x_n)>0
 \end{array} \right.
 \right\rbrace,
 \end{equation}
 \begin{equation}
  U_{I,F,G}^* = \left\lbrace \delta^* \in \sper^*(B_T) \ \left|
 \begin{array}{l}
 \forall j \in F\exists c\in R_{>0}, |y_j|_{\delta^*}>_{\delta^*}c\\
  \exists q \in G,N \in \N \mbox{ s.t. } \forall j \in I,\
 N\nu_{\delta^*}(y_q)> \nu_{\delta*}(y_j) \\
 \forall j \in \{1,\ldots,p\}, \ \forall t \in \{p+1,\ldots,n\},
 \forall N' \in \N, \\\hspace{2cm} \ N'\nu_{\delta^*}(y_j) <\nu_{\delta^*}(y_t)\\
 \nu_{\delta^*}(y_j)>0\ \forall j\in I\cup G
 \end{array} \right.
 \right\rbrace,
\end{equation}
 \begin{equation}
 U_{a,I,F,G} = \left\lbrace \delta \in U_{I,F,G} \ \left|\
 (\nu_\delta(x_1),\ldots,\nu_\delta(x_n)) \underset{\circ}{\sim} (a_1,\ldots,a_n)
\right.\right\rbrace,
 \end{equation}
\begin{equation}
  U_{a,I,F,G}^* = \left\lbrace \delta^* \in  U_{I,F,G}^* \ \left|\
(\nu_{\delta^*}(y_{p+1}),\ldots,\nu_{\delta^*}(y_n))
   \underset{\circ}{\sim} (a_{p+1},\ldots,a_n) \right.\right\rbrace,
 \end{equation}
\begin{equation}
 U_{H,T} = \left\lbrace \delta \in \sper(A) \ \left|
 \begin{array}{l}
 \exists c \in R,\ |x_j|_\delta <_\delta c, \forall j \in H \\
 \exists \epsilon \in R_{> 0}, \ |x_j|_\delta >_\delta \epsilon, \forall j \in T \\
 \nu_\delta(x_1)=\ldots=\nu_\delta(x_p)=0\\
 \nu_\delta(x_{p+1})>0,\ldots,\nu_\delta(x_n)>0
 \end{array} \right.
 \right\rbrace
 \end{equation}
  and
  \begin{equation}
  U_{H,T}^* = \left\lbrace \delta^* \in \sper^*(B) \ \left|
 \begin{array}{l}
\exists q \in T,N \in \N \mbox{ s.t. } \forall j \in H,\
 N\nu_{\delta^*}(y_q)>\nu_{\delta*}(y_j) \\
 \forall j \in \{1,\ldots,p\}, \ \forall t \in \{p+1,\ldots,n\},
 \forall N' \in \N, \\\hspace{2cm} \ N'\nu_{\delta^*}(y_j) <
 \nu_{\delta^*}(y_t)
 \end{array} \right.
 \right\rbrace.
\end{equation}
 We view $U_{I,F,G}$, $U_{a,I,F,G}$ and $U_{H,T}$ (resp. $U_{I,F,G}^*$, $U^*_{a,I,F,G}$ and $U^*_{H,T}$) as topological
 subspaces of $\sper(A)$ (resp. $\sper^*(B_T)$) with the spectral topology. Clearly, for each $I$ and $G$ satisfying
 (\ref{eq:decomp}) we have
 $$
 U_{I,F,G}=\bigcup\limits_{\begin{array}{c}a\in\Gamma^n_+\\a_1=\dots=a_p=0\\a_{p+1}>0,\dots,a_n>0\end{array}}U_{a,I,F,G}
 $$
and
$$
 U^*_{I,F,G}=\bigcup\limits_{\begin{array}{c}a\in\Gamma^n_+\\a_1=\dots=a_p=0\\a_{p+1}>0,\dots,a_n>0\end{array}}
 U^*_{a,I,F,G}
 $$
 Also, we have
 \begin{equation}
 U_{H,T}=\coprod\limits_{\begin{array}{c}\{1,\dots,p\}=I\coprod F\coprod G\\I\subset H,G\subset T\end{array}}U_{I,F,G},
 \label{eq:HTIFG}
 \end{equation}
 and
\begin{equation}
 U^*_{H,T}=\coprod\limits_{\begin{array}{c}\{1,\dots,p\}=I\coprod F\coprod G\\I\subset H,G\subset T\end{array}}U^*_{I,F,G}.
 \label{eq:HTIFG*}
 \end{equation}
  \begin{thm}\label{homeo}
 The map $\psi$ which sends $\delta$ to $\delta^*$, defined above, induces homeomorphisms
\begin{equation}
U_{I,F,G}\tilde\to U_{I,F,G}^*\label{eq:IFG}
\end{equation}
and
\begin{equation}
U_{H,T}\tilde\to U_{H,T}^*.\label{eq:HT}
\end{equation}
If, in addition, $a_{p+1},\ldots,a_n$ are scalewise $\Q$-linearly independent, we also have a homeomorphism
\begin{equation}
U_{a,I,F,G}\tilde\to U_{a,I,F,G}^*.\label{eq:aIFG}
\end{equation}
 \end{thm}

 \noi\textbf{Proof:} To show (\ref{eq:IFG}) and (\ref{eq:aIFG}), we have to prove that
 \begin{equation}
 \psi(U_{I,F,G})\subset U_{I,F,G}^*\label{eq:uIGdansu*}.
 \end{equation}
 \begin{equation}
 \psi(U_{a,I,F,G})\subset U_{a,I,F,G}^*\label{eq:uadansu*},
 \end{equation}
  \begin{equation}
 \psi^{-1}(U^*_{I,F,G})\subset U_{I,F,G}\label{eq:u*IGdansu}
 \end{equation}
 and
\begin{equation}
 \psi^{-1}(U^*_{a,I,F,G})\subset U_{a,I,F,G}\label{eq:u*dansua}.
 \end{equation}
 First, take a point $\delta\in U_{I,F,G}$. By definitions, we have $I_\delta=I$, $F_\delta=F$ and $G_\delta=G$.
For all $j \in F=F_\delta$ there exists $c\in R_{>0}$ such that $|y_j|_{\delta^*}>_{\delta^*}c$ by definition of $F_\delta$
and $B_T$. The condition
 \begin{equation}
 \forall j \in \{1,\ldots,p\}, \ \forall t \in \{p+1,\ldots,n\},
 \forall N' \in \N,  \ N'\nu_{\delta^*}(y_j) <
 \nu_{\delta^*}(y_t)\label{eq:isolated}
 \end{equation}
 is nothing but Proposition \ref{propdelta} (6).

 By Proposition \ref{propdelta} (3), there exist $q \in G$ and $N
 \in \N,\ N > 0$ such that for all $j \in I$ we have
 \begin{equation}
 N\nu_{\delta^*}(y_q)> \nu_{\delta*}(y_j).\label{eq:I<G}
 \end{equation}
 By Proposition \ref{propdelta} (2), we have $\nu_{\delta^*}(y_j)>0$ for all $j\in I\cup G$.
 This completes the proof of the inclusion (\ref{eq:uIGdansu*}).

 Next, assume that $\delta\in U_{a,I,F,G}$ and that $a_{p+1}$, \dots, $a_n$ are $\mathbb Q$-linearly independent. The
 isomorphism
 $$
 (\nu_{\delta^*}(y_{p+1}),\ldots,\nu_{\delta^*}(y_n))
   \underset{\circ}{\sim} (a_{p+1},\ldots,a_n)
   $$
   is given by Proposition \ref{propdelta} (7). This proves the inclusion (\ref{eq:uadansu*}).

   To prove the opposite inclusions, take any $\delta^*\in U_{I,F,G}^*$.
 The existence of $c,N \in R_{>0}$ such that $|x_j|_\delta <_\delta c$ for all $j \in I$
 and
 \begin{equation}
\ |x_j|_\delta >_\delta N,\quad\text{for all }j \in G\label{eq:xGminor}
 \end{equation}
 follow from the facts that $\delta^*$ is bounded,
 $x_j=y_j$ for $j\in I$ and $x_j=1/y_j$ for $j\in G$. For $j\in F$ we have either $x_j=y_j$ or $x_j=\frac1{y_j}$, but
 in both cases the fact that $\delta^*\in U_{I,F,G}^*$ implies the existence of $c_1,c_2\in R_{>0}$ such that
\begin{equation}
 c_1<_\delta|x_j|_\delta<_\delta c_2.\label{eq:c1c2}
\end{equation}
To prove the inclusion (\ref{eq:u*IGdansu}), it remains to prove that
\begin{equation}
\nu_\delta(x_1)=\ldots=\nu_\delta(x_p)=0\label{eq:xp0}
\end{equation}
and
 \begin{equation}
 \nu_\delta(x_t)>0\text{ for all }t\in\{p+1,\ldots,n\}.\label{eq:xGpositif}
 \end{equation}
Equation (\ref{eq:xp0}) is equivalent to saying that
\begin{equation}
 1/|x_j|_\delta\in R_\delta\qquad\text{for }1\le j\le p.\label{eq:inRdelta}
 \end{equation}
 First, if $j\in G$, $|x_j|_\delta=\frac1{|y_j|_{\delta^*}}$ is bounded below by a
 positive constant by (\ref{eq:xGminor}), hence (\ref{eq:inRdelta}) holds for $j\in G$.

If $j\in I$, the assumed existence of $q \in G$ and a positive $N\in \N$ such that for all
 $j\in I$ we have $N\nu_{\delta^*}(y_q) > \nu_{\delta*}(y_j)$
  implies that $|y_j|_{\delta^*}>_{\delta^*}|y_q|_{\delta^*}^N$ by equation (\ref{eq:val2}), in
 other words, $|x_j|_\delta>_\delta1/|x_q|_\delta^N$ or, equivalently, $1/|x_j|_\delta<_\delta |x_q|_\delta^N$. This proves
 (\ref{eq:inRdelta}) for $j\in I$. For $j\in F$, (\ref{eq:inRdelta}) follows from (\ref{eq:c1c2}). Thus
(\ref{eq:inRdelta}) holds for all $j\in\{1,\dots,p\}$, which proves (\ref{eq:xp0}).

 Take an index $t\in\{p+1,\dots,n\}$. To prove (\ref{eq:xGpositif}), it suffices to show that
 \begin{equation}
 1/x_t\notin R_\delta,\label{eq:1/x}
 \end{equation}
 that is, that $1/|x_t|_\delta$ is not bounded above (with respect to
 $\le_\delta$) by any polynomial in $x_1,\ldots,x_n$. By the triangle
 inequality, this is equivalent to saying that $1/|x_t|_\delta$ is not bounded
 above by any monomial in $x_1,\ldots,x_n$, or, equivalently, by any
 element of the form $cx_j^N$ with $j\in\{1,\dots,n\}$, $N\in\N$ and $c\in R$. We prove this last
 statement by contradiction. Suppose there was an inequality of the
 form
 \begin{equation}
 1/|x_t|_\delta<_\delta cx_j^N\label{eq:majoration}
 \end{equation}
 with $N\in\N$, $c\in R$ and $j\in\{1,\ldots,n\}$. Since
 $\nu_{\delta^*}(y_t)>0$, we have $|y_t|_{\delta^*}<_{\delta^*}\epsilon$ for all
 positive $\epsilon\in R$, so $|x_t|_\delta<_\delta\epsilon$ and $1/|x_t|_\delta>1/\epsilon$ for all
 positive $\epsilon\in R$. On the other hand, if $j\in
 I\cup\{p+1,\ldots,n\}$, we have $\nu_{\delta^*}(y_j)>0$, hence
 $|x_j|_\delta=|y_j|_{\delta^*}<_{\delta^*}\theta$ for all positive $\theta\in R$ and if $j\in F$ then $|x_j|_\delta$ is
 bounded above by a constant from $R$ by (\ref{eq:c1c2}). This proves
 that $j\in G$ in (\ref{eq:majoration}).

 Now, the hypotheses (\ref{eq:isolated}) implies that for any constant
 $d\in R$ and any $N'\in\N$ we have $d|y_j|_{\delta^*}^{N'}>_{\delta^*}|y_t|_{\delta^*}$, so
 $d/|x_j|_\delta^{N'}>_\delta|x_t|_\delta$, which contradicts (\ref{eq:majoration}). This
 completes the proof of (\ref{eq:1/x}) and (\ref{eq:xGpositif}). The inclusion (\ref{eq:u*IGdansu}) is proved.

Assume that $\delta^*\in U_{a,I,G}^*$. To prove the inclusion (\ref{eq:u*dansua}), it remains to prove the isomorphism
 \begin{equation}
 (\nu_\delta(x_1),\ldots,\nu_\delta(x_n)) \underset{\circ}{\sim}
 (a_1,\ldots,a_n).\label{eq:ogm}
 \end{equation}
By Proposition \ref{propdelta} (5), (\ref{eq:xGpositif}), the
 assumed scalewise $\Q$-linear independence of $a_{p+1},\ldots$, $a_n$ and
 the Remark following Definition \ref{scalewise},
 $\nu_\delta(x_{p+1}),\ldots,\nu_\delta(x_n)$ are also scalewise
 $\Q$-linearly independent. Now (\ref{eq:ogm}) follows from Proposition \ref{propdelta} (7). The inclusion
 (\ref{eq:u*dansua}) is proved.

 Finally, the homeomorphism (\ref{eq:HT}) follows from (\ref{eq:IFG}), (\ref{eq:HTIFG}) and (\ref{eq:HTIFG*}) by letting
 $I,F,G$ run over all the triples of disjoint subsets, satisfying (\ref{eq:decomp}), (\ref{eq:IdansH}) and
 (\ref{eq:GdansT}), such that $I$ is empty whenever $G$ is empty.
 \hfill $\Box$
 \medskip

 Of course, Theorem \ref{homeo} is true with $\{1,\dots,p\}$ replaced by any other subset of $\{1,\dots,n\}$. In the
 next Corollary we drop the assumption $(\ref{eq:decomp})$ and let $I,F,G$ run over all the triples of disjoint subsets
 of $\{1,\dots,n\}$ such that $I=\emptyset$ whenever $G=\emptyset$. Similarly, $H,T$ will run over all the pairs of
disjoint subsets of $\{1,\dots,n\}$.
\begin{cor} We have finite coverings
$$
\sper\ A=\coprod\limits_{I,F,G}U_{I,F,G},
$$
and
$$
\sper\ A=\bigcup\limits_{H,T}U_{H,T}
$$
For each $I,F,G$ as above, the set
$U_{I,F,G}$ is homeomorphic to the subset $U^*_{I,F,G}$ of the set $\sper^*B_G$ of finite points of
$\sper\ B_G$. For each $H,T$ as above, the set
$U_{H,T}$ is homeomorphic to the subset $U^*_{H,T}$ of the set $\sper^*B_T$ of finite points of $\sper\
B_T$.
\end{cor}
\begin{rek} The assumption of scalewise $\mathbb Q$-linear independence of $a_{p+1},\dots,a_n$,
is needed in Theorem \ref{homeo} only for the inclusion (\ref{eq:u*dansua}). The usual $\mathbb Q$-linear
independence is needed for the inclusion (\ref{eq:uadansu*}) and for Proposition \ref{propdelta}. Although at first glance
these assumptions seem rather restrictive, we remark that any point $\delta\in\sper\ A$ can be transformed into one
for which these assumptions hold by a sequence of blowings up. We refer the reader to Corollary 6.2 of \cite{LMSS} for
details. Corollary 6.2 of \cite{LMSS} shows how to achieve usual $\mathbb Q$-linear independence of $a_{p+1},\dots,a_n$,
but it also works for scalewise $\mathbb Q$-linear independence after some minor and obvious modifications.
\end{rek}
 \noi\textbf{Example.} Let $n=5$. Let $\delta\in\sper\ A$ be the point given by the following semi-curvette. We let
 $\Gamma=\mathbb Z^2_{\text{lex}}$ and $k_\delta=R(z,w)$, where $z$ and $w$ are independent variables. Let the
 total order on $k_\delta$ be given by the following inequalities:
\begin{eqnarray}
0&<_\delta&w<_\delta c<_\delta z\quad\text{ for all }\ c\in R_{>0}\\
\frac1{w^N}&<_\delta&z\qquad\qquad\qquad\text{for all }N\in\mathbb N.
\end{eqnarray}
As usual, we define the total order on $k_\delta\left(\left(t^\Gamma\right)\right)$ by declaring $t$ to be positive. Define
the map $\delta:A\rightarrow k_\delta\left(\left(t^\Gamma\right)\right)$ by
 \begin{eqnarray}
\delta(x_1)&=&w\\
 \delta(x_2)&=&1+t^{(0,1)}\\
 \delta(x_3)&=&z\\
 \delta(x_4)&=&t^{(1,0)}\\
 \delta(x_5)&=&zt^{(1,0)}.
 \end{eqnarray}
We have $\nu_\delta(x_1)=\nu_\delta(x_2)=\nu_\delta(x_3)=0$,
\begin{equation}
\nu_\delta(x_4)=\nu_\delta(x_5)=(1,0)>0,\label{eq:345}
\end{equation}
so $p=3$. Moreover,
$I_\delta=\{1\}$, $F_\delta=\{2\}$, $G_\delta=\{3\}$. Let $T=G_\delta$ and let $\delta^*=\psi(\delta)\in\sper^*B_T$. We
have $\Gamma_{\delta^*}=\mathbb Z^4_{\text{lex}}$ and $k_{\delta^*}=R$.
The semi-curvette $\delta^*$ can be defined by the map
\begin{eqnarray}
 \delta^*(y_1)&=&t^{(0,0,0,1)}\\
 \delta^*(y_2)&=&1+t^{(0,1,0,0)}\\
 \delta^*(y_3)&=&t^{(0,0,1,0)}\\
 \delta^*(y_4)&=&t^{(1,0,0,0)}\\
 \delta^*(y_5)&=&t^{(1,0,1,0)}.
 \end{eqnarray}
 In this example, $\nu_\delta(x_4)$ and $\nu_\delta(x_5)$ are not $\Q$-linearly independent
 (\ref{eq:345}), and the conclusion of Proposition \ref{propdelta} does not hold: we do \textit{not} have the equivalence
 $$
 (\nu_{\delta^*}(y_4),\nu_{\delta^*}(y_5)) \underset{\circ}{\sim}
  (\nu_{\delta}(x_4),\nu_{\delta}(x_5)).
 $$
Let $A'=R[x'_1,x'_2,x'_3,x'_4,x'_5]$. Consider the map $\pi:A\rightarrow A'$ defined by
\begin{eqnarray}
 \pi(x_j)&=&x'_j\quad\text{ for }j\in\{1,2,3,4\},\\
 \pi(x_5)&=&x'_4x'_5.
\end{eqnarray}
 Let $\delta'$ be the unique preimage of $\delta$ under the
 natural map $\pi^*:\sper\ A'\rightarrow\sper\ A$ of the real spectra, induced by $\pi$ (in the terminology of \cite{LMSS},
 $\pi$  is an affine monomial blowing up along the ideal $(x_4,x_5)$ with respect to $\delta$ and $\delta'$ is the
 transform of $\delta$ by $\pi$). Explicitly, we have $\Gamma_{\delta'}=\Z^2_{\text{lex}}$, $k_{\delta'}=R(z,w)$, as above,
 and  $\delta$ is given by the semi-curvette
\begin{eqnarray}
\delta(x_1)&=&w\\
 \delta(x_2)&=&1+t^{(0,1)}\\
 \delta(x_3)&=&z\\
 \delta(x_4)&=&t^{(1,0)}\\
 \delta(x_5)&=&z.
 \end{eqnarray}
 This is an example of the fact that every point $\delta\in\sper\ A$ can be transformed, after a sequence
$\sper\ A'\rightarrow\sper\ A$ of affine monomial blowings up with respect to $\delta$, into a point $\delta'\in\sper\ A'$
such that the non-zero elements of the set $\{\nu_{\delta'}(x_1),\dots,\nu_{\delta'}(x_n)\}$ are (scalewise) $\mathbb
Q$-linearly independent.


\begin{thebibliography}{99}
\bibitem{Baer} R. Baer, Uber nicht-archimedisch geordnete K\"{o}rper (Beitrage zur Algebra). Sitz. Ber. Der Heidelberger
Akademie, 8 Abhandl. (1927).
\bibitem{BCR} J. Bochnak, M. Coste, M.-F. Roy, {\em G\'eom\'etrie
    alg\'ebrique r\'eelle.} Springer--Verlag, Berlin 1987.
\bibitem{Fuc} L. Fuchs, {\em Telweise geordnete algebraische Strukturen.} Vandenhoeck and Ruprecht, 1966.
\bibitem{Kap} I. Kaplansky, {\em Maximal fields with valuations I.}
  Duke Math. J., 9:303--321 (1942).
\bibitem{Kap2} I. Kaplansky, {\em Maximal fields with valuations II.}
  Duke Math. J., 12:243--248 (1945).
\bibitem{Krull} W. Krull, Allgemeine Bewertungstheorie, J. Reine Angew. Math. 167, 160--196 (1932).
\bibitem{LMSS} F. Lucas, J. J. Madden, D. Schaub, M. Spivakovsky {\em On connectedness of sets in the real spectra of
polynomial rings}, arXiv AG/0601671
\bibitem{Pre} A. Prestel {\em Lectures on formally real fields}, Lecture Notes in Math., Spriniger--Verlag---Berlin,
Heidelberg, New York, 1984.
\bibitem{PC} S. Priess-Crampe, {\em Angeordnete strukturen: gruppen, orper, projektive Ebenen}, Springer--Verlag---Berlin,
Heidelberg, New York, 1983.
\bibitem{Spi2} M. Spivakovsky, {\em A solution to Hironaka's polyhedra
    game.} Arithmetic and Geometry, Vol II, Papers dedicated to I. R.
  Shafarevich on the occasion of his sixtieth birthday, M. Artin and
  J. Tate, editors, Birkh\"{a}user, 1983, pp. 419--432.
 \bibitem{ZS} O. Zariski, P. Samuel, Springer--Verlag---Berlin,
Heidelberg, New York, 1960.
\end{thebibliography}
\end{document}